\documentclass[final,1p,times]{elsarticle}
\usepackage{siunitx}
\usepackage{amssymb}
\usepackage{threeparttable}
\usepackage{booktabs}
\usepackage{amsmath}
\usepackage{hyperref}
\usepackage{graphicx}
\usepackage{algorithm}
\usepackage{algpseudocode}
\usepackage{mathrsfs}
\usepackage{lineno}
\usepackage{amsthm}
\journal{Journal of Computational Science}

\begin{document}

\begin{frontmatter}



\title{SRL-Assisted AFM: Generating Planar Unstructured Quadrilateral Meshes with Supervised and Reinforcement Learning-Assisted Advancing Front Method}


\author[label1]{Hua Tong}
\author[label1]{Kuanren Qian}
\author[label1,label2,label3]{Eni Halilaj}
\author[label1,label2]{Yongjie Jessica Zhang}

\address[label1]{Department of Mechanical Engineering, Carnegie Mellon University, 5000 Forbes Ave, Pittsburgh, PA 15213, USA}
\address[label2]{Department of Biomedical Engineering, Carnegie Mellon University, 5000 Forbes Ave, Pittsburgh, PA 15213, USA}
\address[label3]{Robotics Institute, Carnegie Mellon University, 5000 Forbes Ave, Pittsburgh, PA 15213, USA}

% \affiliation[label1]{organization={Department of Mechanical Engineering, Carnegie Mellon University},%Department and Organization
%             addressline={5000 Forbes Ave},
%             city={Pittsburgh},
%             state={PA 15213},
%             % postcode={15213},
%             country = {USA}}
% \affiliation[label2]{organization={Department of Biomedical  Engineering, Carnegie Mellon University},%Department and Organization
%             addressline={5000 Forbes Ave},
%             city={Pittsburgh},
%             state={PA 15213},
%             % postcode={15213},
%             country = {USA}}
% \affiliation[label3]{organization={Robotics Institute, Carnegie Mellon University},%Department and Organization
%             addressline={5000 Forbes Ave},
%             city={Pittsburgh},
%             state={PA 15213},
%             % postcode={15213},
%             country = {USA}}

\begin{abstract}
High-quality mesh generation is the foundation of accurate finite element analysis. Due to the vast interior vertices search space and complex initial boundaries, mesh generation for complicated domains requires substantial manual processing and has long been considered the most challenging and time-consuming bottleneck of the entire modeling and analysis process. In this paper, we present a novel computational framework named ``SRL-assisted AFM" for meshing planar geometries by combining the advancing front method with neural networks that select reference vertices and update the front boundary using ``policy networks." These deep neural networks are trained using a unique pipeline that combines supervised learning with reinforcement learning to iteratively improve mesh quality. First, we generate different initial boundaries by randomly sampling points in a square domain and connecting them sequentially. These boundaries are used for obtaining input meshes and extracting training datasets in the supervised learning module. We then iteratively improve the reinforcement learning model performance with reward functions designed for special requirements, such as improving the mesh quality and controlling the number and distribution of extraordinary points. Our proposed supervised learning neural networks achieve an accuracy higher than $98\%$ on predicting commercial software. The final reinforcement learning neural networks automatically generate high-quality quadrilateral meshes for complex planar domains with sharp features and boundary layers.
\end{abstract}




\end{frontmatter}

\section{Introduction}
As outlined in NASA's Vision 2030, mesh generation constitutes one of the six crucial research directions and holds significance in numerical simulations \cite{slotnick2014cfd}. However, mesh generation remains a major bottleneck due to algorithmic complexity, poor error estimation capabilities, and intricate geometries \cite{zhang2018geometric, zhang2013challenges}. Quadrilateral (Quad) meshes are typically chosen over triangular meshes in applications such as texturing, simulation using finite elements, and B-spline fitting due to their attractive tensor-product nature and smooth surface approximation. Quad mesh generation has been a significant research topic for decades. Yet, existing quad mesh generation techniques rely heavily on pre-processing or post-processing to maintain good mesh quality and require heuristic expertise in algorithm construction. Pre-processing involves creating optimal vertex locations \cite{remacle2013frontal} and breaking down complex domains into regular components \cite{liu2017distributed}. Post-processing is performed to clean inverted or irregularly connected elements. The operations include mesh adaptation \cite{verma2015robust}, splitting, swapping, and collapsing elements as iterative topological alterations \cite{docampo2019towards}, as well as singularity reduction \cite{verma2016alphamst}. However,  these additional mesh quality improvement procedures are computationally complex and inefficient.

There are two types of quad meshes: structured \cite{thompson1998handbook} and unstructured \cite{owen1998survey, teng2000unstructured}. All elements in a structured mesh are arranged in a regular pattern before being mapped to the user-defined boundaries. The resulting mesh quality may be poor for complex boundaries. In unstructured meshes, interior node valence numbers are relaxed, allowing for greater flexibility in mesh construction. There are two methods for generating unstructured quad meshes: indirect and direct. The indirect method first triangulates the domain, then employs edge midpoints and face centers to transform each triangular element into four quads to optimize triangular and quad surface mesh quality \cite{garimella2004triangular}. An advancing front method (AFM) was proposed to generate an all-quad mesh from triangles \cite{lee1994new}. The initial front of the mesh is defined by delineating the triangle edges at the boundary. A sequence of paired triangles is systematically merged along the front, progressively moving toward the interior. Q-Morph is another AFM-based method that can effectively decrease the number of irregular nodes by performing local edge swapping and inserting additional nodes \cite{owen1999q}. However, indirect methods require an intermediate triangular mesh, which is prone to instability and has a restricted number of vertices. The direct method, on the other hand, bypasses triangulation entirely and instead generates quad elements directly, avoiding those potential problems. 

Many direct methods have been proposed in recent years. The paving method initiates from the boundary and proceeds inward by arranging complete rows of elements on the front boundary once at a time \cite{blacker1991paving}. The quadtree/hexagon-based methods are hierarchical approaches to mesh generation subdividing a region into quad \cite{baehmann1987robust, liang2009guaranteed} or hexagonal \cite{liang2011hexagon} cells recursively based on geometric criteria. In \cite{liang2012matching}, an individual interior or exterior mesh is generated at a time and matched at the shared boundary. Mesh quality can be improved via face swapping, edge removal, and geometric flow-based smoothing \cite{zhang2009surface}. Octree-based iso-contouring methods analyze each interior grid vertex and generate a dual mesh of the background grids for any complicated single-material and multiple-material domains \cite{zhang20053d, zhang2006adaptive, zhang2010automatic}. In biomedical applications, vascular blood flow simulation using isogeometric analysis (IGA) needs high-order elements like T-splines \cite{wei2022analysis}, subdivisions \cite{wei2015truncated, wei2021tuned}, and THB-splines \cite{wei2017truncated}. All prefer high-quality unstructured quad and hexahedral meshes because they can be converted into standard or rational T-splines, which possess $C^2$-continuity across the entirety of the surface, excluding local regions proximal to the extraordinary points (EPs) \cite{wang2011converting, wang2012converting}. Together with local refinement in IGA, the computational efficiency can be greatly improved.

Our research is based on the AFM, one of the most widely-used direct methods \cite{docampo2019towards}. AFM is a greedy algorithm that iteratively creates mesh nodes from the input boundaries to the interior. Each iterative process consists of three steps: (1) selecting a line segment from the front set that separates the meshed domain and the unmeshed domain; (2) connecting a new mesh node or existing mesh nodes to the base segment to generate a high-quality quad element; and (3) updating the front set until the entire domain is meshed. AFM can produce a high-quality mesh, but it is inefficient since it requires numerous intersection calculations \cite{lohner1988generation, guo2021improved}. Several literatures attempt to integrate mesh generation with machine learning (ML) modules to create new meshing algorithms. Training a reinforcement learning neural network on an input boundary multiple episodes to generate a final good-quality mesh \cite{pan2023reinforcement} is one example. Deep learning is also used to learn the progress direction and step size of triangle mesh generation \cite{lu2022new}. However, these methods require intersection detection in the middle stage or after mesh generation, which is computationally expensive. In addition, these methods were only tested on simple boundaries and did not provide open-source codes to test their generalization for complex domains. 

To tackle complex boundaries and satisfy special requirements in generating quad meshes, we propose a new computational framework named ``SRL-assisted AFM" for quad mesh generation using AFM assisted by supervised learning (SL) and reinforcement learning (RL). We train the SL module with the dataset extracted  from input meshes. The RL module automatically generates high-quality meshes with designed reward functions for various complex geometries and back-propagates neural network weights based on high-quality training datasets. The proposed SRL-assisted AFM framework is capable of meshing complex new boundaries efficiently. Users can also add boundary layers and optimize the number and distribution of EPs in the mesh. The main contributions of this paper include:
\begin{itemize}
\item Integrating AFM, SL, and RL into a new computational framework to generate high-quality quad meshes for planar domains with complex boundaries;
\item Eliminating the need for quality improvement and intersection detection during and after mesh generation; and
\item Preserving high Jacobian, low aspect ratio, low EP number with sharp features, unbalanced seeds, conformal boundaries, and boundary layers.
\end{itemize}

The remainder of this paper is structured as follows. In Section 2, we overview the comprehensive framework of AFM mesh generation as an SL-RL problem. In Section 3, we discuss the detailed AFM, SL, and RL modules, along with the comprising action, state, and reward settings. In Section 4, we present numerical results and discuss our findings. We conclude by summarizing our contributions and proposing future work in Section 5.

\section{Overview of the SRL-Assisted AFM Framework}
\begin{figure}[ht]
\centering
\includegraphics[width=\linewidth]{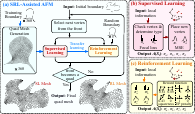}
\vspace{-6mm}
\caption{Neural network pipeline and architecture. (a) The SRL-assisted AFM framework utilizes SL and RL neural networks to generate new vertices following the AFM scheme. (b) SL neural networks learn from ANSYS-generated quad meshes and generate a new quad on the front. (c) RL neural networks use squareness and EP penalty as the reward functions and are trained on the evolving front to maximize the reward feedback.}
\label{fig:bigpicture}
\end{figure}

We combine the AFM with the SL and RL modules to automatically generate high-quality quad meshes (Figure~\ref{fig:bigpicture}). AFM generates one quad element at a time based on the front boundary information and evolves the front at each time step. The meshing process is complete when the evolving front becomes a quad element. In the literature, rule-based algorithms are adopted at each step to generate a new element \cite{lohner1988generation, seveno1997towards}. Here we replace rule-based algorithms with policy neural networks (SL and RL), which is capable of approximating any complex functions \cite{hornik1989multilayer}.

Our SRL-assisted AFM framework (Figure \ref{fig:bigpicture}(a)) includes a training procedure (blue arrows) and a testing procedure (black arrows). In the training procedure, we first randomly generate 360 planar training boundaries by connecting randomly placed vertices in a square domain. Then, we use commercial software ANSYS to generate quad meshes based on these boundaries and collect $3.5M$ quad elements for our training dataset. The four SL policy neural networks (Figure \ref{fig:bigpicture}(b)) $\left\{\pi_a, \pi_b, \pi_c, \pi_d\right\}$ take local information around a selected vertex as input and work together sequentially to update the front: (1) collecting local information from the vertex with the smallest angle on the front, (2) sending local information to $\pi_a$ and selecting the reference point, (3) sending reference point local information to $\pi_b$ and getting the updating type, and (4) sending reference point local information to $\pi_c$ and $\pi_d$ to generate new interior vertices if the $\pi_b$ updating type requires inserting new points. We use input meshes results to obtain the optimized policy network weights via AFM and SL to circumvent random weights-induced poor performance. 

To further improve the framework, we combine AFM with RL neural networks (Figure \ref{fig:bigpicture}(c)). We consider the AFM as a Partially Observable Markov Decision Process \cite{cassandra1998survey} and combine it with RL. We transfer trained SL neural networks to RL neural networks with the same architecture denoted as $\left\{\pi'_a, \pi'_b, \pi'_c, \pi'_d\right\}$ and train them on $5$ randomly selected initial boundaries (Figure \ref{fig:bigpicture}(a)). At each time step $t$ during the RL training process, the environment (current front) sends a state $S_t$ to the neural network. Then the neural network samples an updating action $A(S_t)$ and receives a reward (mesh quality metrics, Figure \ref{fig:bigpicture}(c)) as feedback of its action. In our implementation, the reward function measures the squareness of elements as well as the number and distribution of EPs. The higher the squareness reward, the closer a quad element is to a square. For the EP penalty, all regular vertices receive reward ``$1$," while EPs receive penalty ``$0$." A pair of EPs that are adjacent to each other are denoted as close EPs or cEPs. All cEPs also receive the ``$0$" penalty. These three steps (receiving the local information, sampling actions with neural network predictions, and updating the front) are iterated until the evolving front becomes a quad, signaling the meshing process is completed. As we mesh the same domain multiple times, the mesh (dataset) quality improves, and we update the neural network weights with the latest high-quality dataset. In the end, we obtain RL neural networks and meshes superior over input meshes in various mesh quality metrics.

\section{Methodology}
Our SRL-assisted AFM pipeline combines three fundamental modules (AFM, SL and RL) together. The technical details of each module are explained below.

\subsection{SRL-Assisted AFM}
\begin{figure}[ht]
\centering
\includegraphics[width=\linewidth]{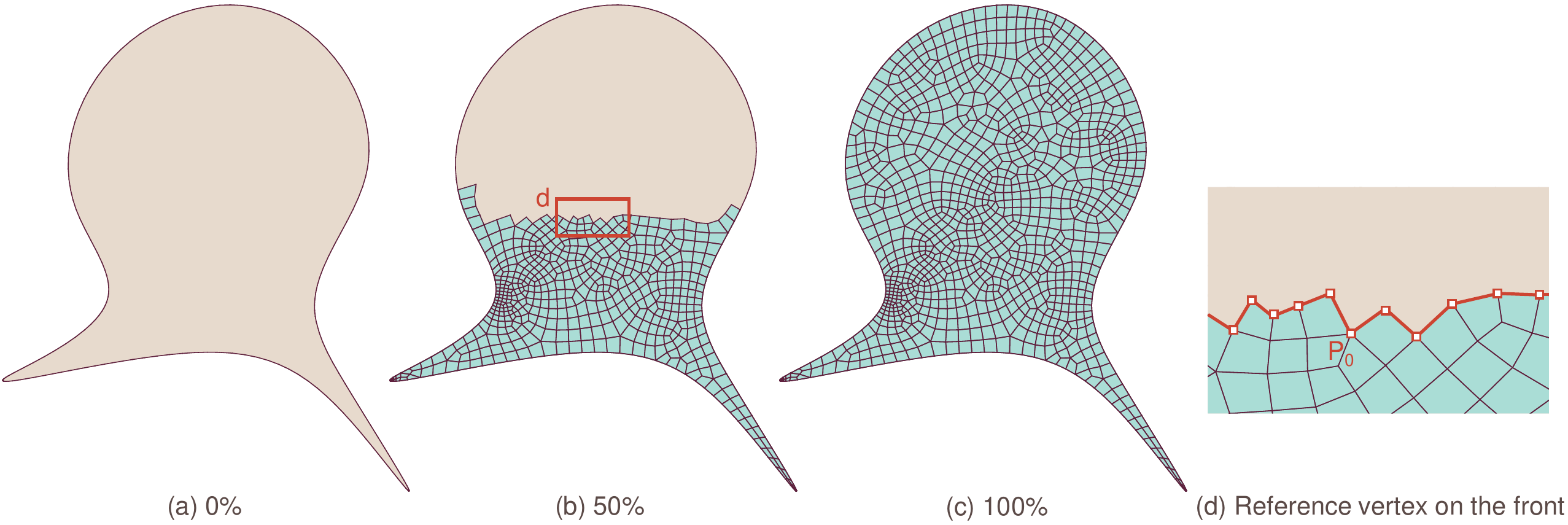}
\vspace{-6mm}
\caption{Meshing a planar domain with advancing front method. (a) The initial domain boundary. (b) The meshing process at $50\%$ completeness. (c) The final mesh. (d) A zoom-in picture of the reference vertex ($P_0$) and the front (red line) at $50\%$ completeness.}
\label{fig:afminverse}
\end{figure}

\begin{algorithm}[ht]
    \textbf{Input}: Planar closed manifold boundaries $B^0_i, i=1, 2, \cdots, L$\\
    \textbf{Output}: Interior or exterior quad meshes \\
    \textbf{Goal}: Fill the interior or exterior domain enclosed by $B^0_i, i=1, 2, \cdots, L$ with all quad elements
    \caption{SRL-assisted AFM}
    \label{alg:AFM}
    \begin{algorithmic}[1]
        \For{boundary number $i = 1, 2, \cdots, L$}
        \State Assign adaptive seeds to $B^0_i$
        \For{boundary number $j = 1, 2, \cdots, L, j\neq i$}
        \If {$B^0_i$ and $B^0_j$ are conformal boundaries}
        \State Calculate adaptive mesh coefficients on both boundaries (denoted as $k_i$ and $k_j$)
        \State Assign new seeds with coefficients $k_{new} = \frac{k_i + k_j}{2}$ to both boundaries
        \EndIf
        \EndFor
        \EndFor
        \For{time step $t = 1, 2, \cdots, M$}
        \For{boundary number $i = 1, 2, \cdots, L$}
        \For{vertex number $j = 1, 2, \cdots, N$}
        \State Save local information around vertex $P_j$ from $B^t_i$ as input tensor $S_j$
        \State Input $S_j$ to $\pi_a$ (or $\pi'_a$) to judge whether $P_j$ can be the reference vertex
        \EndFor
        \State Denote the selected reference vertex and its local information as $P_0, S_t$
        \State Input $S_t$ to $\pi_b$ (or $\pi'_b$) to obtain the updating type $T_t$
        \If {$T_t=1$}
        \State Input $S_t$ to $\pi_c$ (or $\pi'_c$) to obtain one new vertex $P^{new}_1$
        \EndIf
        \If {$T_t=4$}
        \State Input $S_t$ to $\pi_d$ (or $\pi'_d$) to obtain two new vertices $P^{new}_{1, 2}$
        \EndIf
        \State Form a new quad element and update front boundary $B_i^t$
        \EndFor
        \EndFor
        \State Add boundary layers
    \end{algorithmic}
\end{algorithm}

In the paper we employ AFM, which iteratively generates one quad element at a time to fill the entire domain (Figure \ref{fig:afminverse}). Our implementation is shown in Algorithm \ref{alg:AFM}. It takes planar closed manifold boundaries $B^0_i, i=1, 2, \cdots, L$ as input and outputs an all-quad mesh. After specifying the initial boundary, we assign seeds to the boundary based on a size function $s(i)$, which is proportional to the local boundary curvature $\rho(i)$. We use the circumcircle of three local vertices $P_{i-1}, P_i$, and $P_{i+1}$ to estimate $\rho(i)$. Then we obtain $s(i) = k(i) \rho(i)$, and $k(i)$ is calculated by
\begin{equation}
k(i) = \frac{l_{tot}}{s_a\sum_{i=1}^N\frac{l(i)}{\rho(i)}}, \quad \textit{where}\ \rho(i) \neq 0. \label{eqn: coefficient k equation} \\
\end{equation}
$\noindent s_{a}$ is an approximate mesh size (the default $s_a = 0.05$). We calculate the whole front boundary length $l_{tot}=\sum_{i=1}^Nl_i$, where $l_i$ is the length of edge $P_iP_{i+1}$.

The seed assignment procedure supports conformal boundaries. Suppose we have a pair of conformal boundaries $B_1$ and $B_2$. At the same node $i$ on $B_1$ or $j$ on $B_2$, the local front boundary curvature on both boundaries is the same, $\rho_1(i)=\rho_2(j)$, because the three local vertices are on the shared curve. However, the coefficients on $B_1$ and $B_2$ are different or $k_1(i)\neq k_2(j)$ because $B_1$ and $B_2$ may have different $l_{tot}$ and $s_a$ values. We assign the new coefficient $\frac{k_1(i)+k_2(j)}{2}$ to this node on the shared curve. In this way, seeds on both boundaries conform exactly to each other.

\begin{figure}[ht]
\centering
\includegraphics[width=\linewidth]{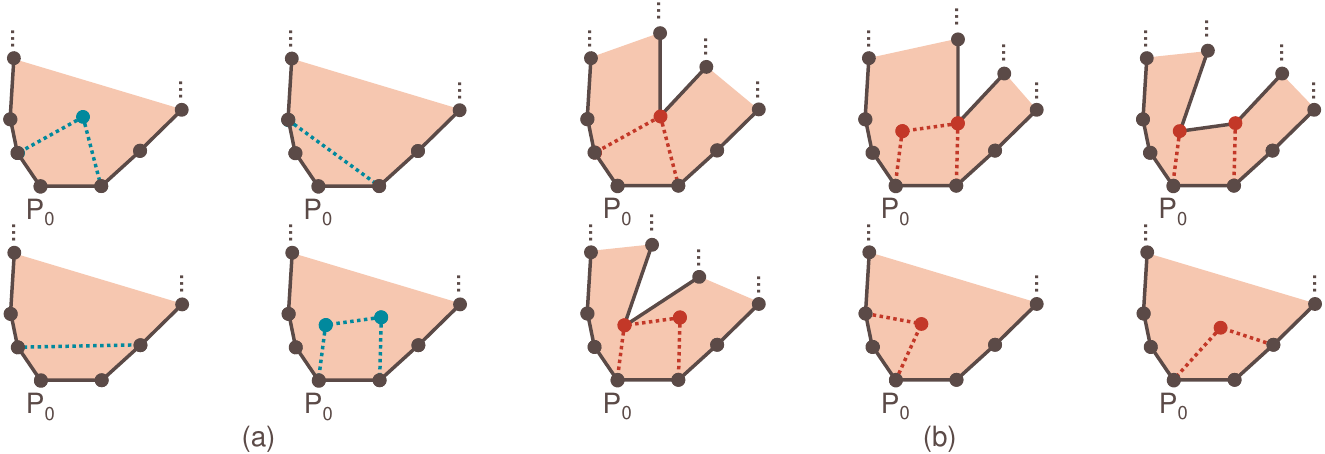}
\vspace{-6mm}
\caption{The judgment principle of the binary classification neural network $\pi_a$ when labeling the ground truth dataset. (a) Four cases (corresponds to the four updating types) when $\pi_a$ accepts $P_0$ to be the reference vertex. (b) Six example cases when $\pi_a$ rejects $P_0$ to be the reference vertex. }
\label{fig:accept_reject_types}
\end{figure}

After assigning the seeds, we begin to mesh the domain. At each time step $t$, we select a reference vertex $P_0$ from the front boundary (the reference vertex position determines where we update the front boundary), collect local information $S_t$ around $P_0$, and take action $A(S_t)$ to update the front. In our SRL-assisted AFM framework, we have four neural networks $\left\{\pi_a, \pi_b, \pi_c, \pi_d\right\}$ in SL and $\left\{\pi'_a, \pi'_b, \pi'_c, \pi'_d\right\}$ in RL that identify the proper $P_0$ and generate a new quad around $P_0$ at each time step. After iterating through all vertices on the front, a vertex is chosen as the reference $P_0$ based on the judgment given by the binary classification neural network ($\pi_a$ in SL, $\pi'_a$ in RL). As shown in Figure \ref{fig:accept_reject_types}(a), $\pi_a$ accepts $P_0$ to be the reference vertex for these four cases, which corresponds to the four updating types. Note that Types $1$ and $4$ insert one and two new vertices, respectively, and Types $2$ and $3$ connect two existing vertices on the front. Figure \ref{fig:accept_reject_types}(b) show six example cases, where $\pi_a$ rejects the point when the formed quad element separates the domain into two subdomains or only one new point is generated and connects with $P_0$ (introducing more complex cases). After this, an updating type is selected based on the four-class classification neural network ($\pi_b$ in SL, $\pi'_b$ in RL). $\pi_c$ and $\pi_d$ ($\pi'_c$ and $\pi'_d$ in RL) are then used to insert one or two vertices in the domain to form a new quad element based on selected types. The meshing process is done only when the number of edges on the front becomes $4$, forming the last quad.

\begin{figure}[ht]
\centering
\includegraphics[width=\linewidth]{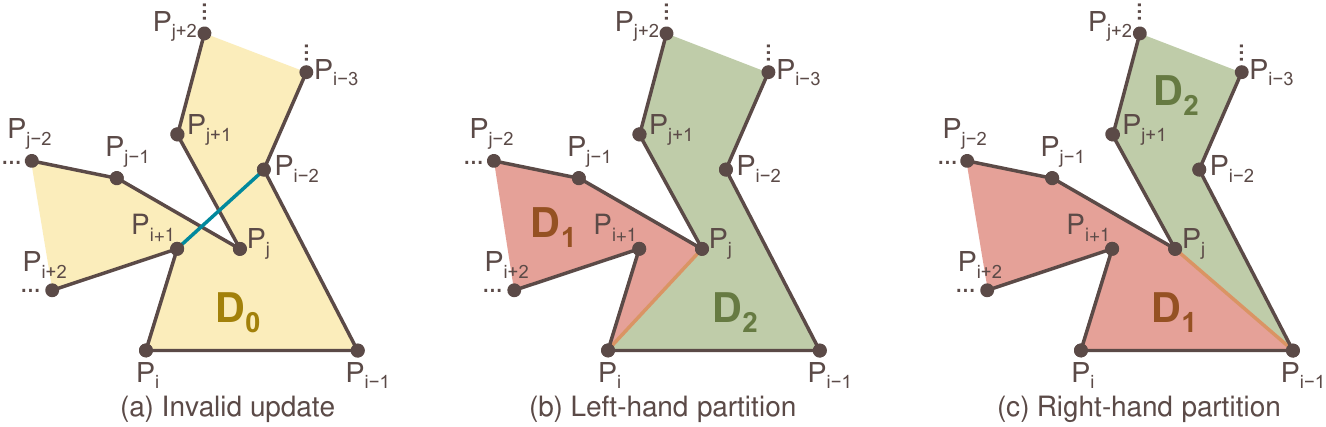}
\vspace{-6mm}
\caption{Two special situations when a partition operation is necessary to split the domain $D_0$. (a) Assuming that $P_i$ is the reference vertex, the next update, which seals the line segment $P_{i+1}P_{i-2}$, is invalid because $P_{i+1}P_{i-2}$ intersects with the remaining front boundary. (b) $D_0$ is partitioned along $P_iP_j$ if the resulting two new subdomains $D_1, D_2$ both have an even number of edges. (c) $D_0$ is partitioned along $P_{i-1}P_j$ if the edge numbers of both new subdomains $D_1, D_2$ are even.}
\label{fig:proof}
\end{figure}

During the mesh generation process, when the updated line segments intersect with the current front, we can correct the error by partitioning the front into two new fronts. Figure \ref{fig:proof} shows a local region of the evolving front, and neural network $\pi_b$ (or $\pi'_b$) selects Type $2$ or Type $3$ classification. Normally, we connect line segment $P_{i+1}P_{i-2}$ (the blue line segment) that corresponds to Type $2$, or $P_{i+2}P_{i-1}$ that corresponds to Type $3$. However, if $P_{i+1}P_{i-2}$ intersects with the remaining front boundary ($P_{j-1}P_j$ and $P_{j+1}P_j$), we partition the front boundary into two new front boundaries which form two subdomains. Eventually, Algorithm \ref{alg:AFM} will fill these two subdomains with quad elements. Theorem \ref{the:the1} discusses special intersection situations when separating the original front boundary into two boundaries is needed to continue the meshing process.

\newtheorem{theorem}{Theorem}
\begin{theorem}
$\forall$ A planar domain $D_0$ with an even number of edges $N_0 \geq 4, \exists$ a partition method to split $D_0$ into two subdomains $D_1, D_2$, both of which have an even number of edges $N_1$ and $N_2$, satisfying $N_1+N_2-2 = N_0$.
\begin{proof}
In Figure \ref{fig:proof}, assume there are multiple line segments $(P_{j-1}P_j$ and $P_jP_{j+1})$ on $D_0$ that intersect with $P_{i-2}P_{i+1}$, we can define a $\text{set}\ \mathcal{P}$ such that $\text{all}\ P_j \in \mathcal{P}$. In set $\mathcal{P}$, we can always find a vertex $P_j$ that is the closest to the line segment $P_iP_{i-1}$:

\begin{align}
P_{j} = \operatorname*{arg min}_{j} dist(P_j, P_iP_{i-1}) = \operatorname*{arg min}_{j} \frac{\left| \overrightarrow{P_iP_{i-1}} \times \overrightarrow{P_iP_j} \right|}{|\overrightarrow{P_iP_{i-1}}|}.
\end{align}
There exist no further line segments on $D_0$ that intersect with $P_iP_j$ or $P_{i-1}P_j$, which means $D_0$ can be partitioned along either $P_iP_j$ or $P_{i-1}P_j$. We partition $D_0$ into two subdomains along $P_iP_j$ or $P_{i-1}P_j$. The operation forms two subdomains with the number of edges $N_1$ and $N_2$, respectively. The choice of partitioning along $P_iP_j$ or $P_{i-1}P_j$ depends on which choice ensures both $N_1$ and $N_2$ are even. Otherwise, $D_1$ and $D_2$ cannot be filled by all quads. After the partition, one new edge is created and shared by $D_1$ and $D_2$. Therefore, we have $N_1+N_2-2=N_0$.
\end{proof}
\label{the:the1}
\end{theorem}

\begin{figure}[ht]
\centering
\includegraphics[width=\linewidth]{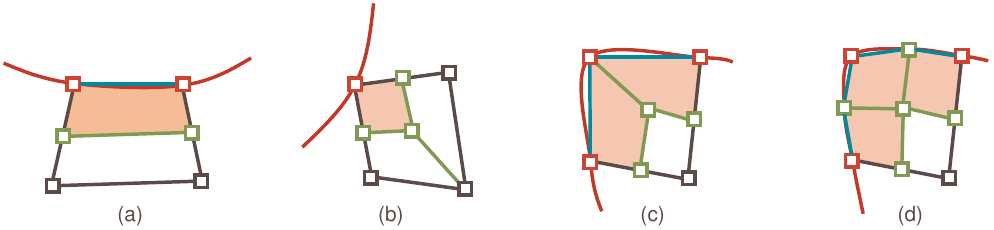}
\vspace{-6mm}
\caption{Four templates of adding the boundary layer. The quad element shares one edge (a), one vertex (b) and two adjacent edges (c, d) with the boundary. A valence-3 EP is introduced in (b, c), while two new vertices are inserted on the boundary in (d).}
\label{fig:boundarylayer}
\end{figure}

After the meshing, we can add boundary layers to the internal boundary for fluid mechanics simulations. We construct the boundary layer by splitting the elements along the boundary. Four templates are implemented (Figure \ref{fig:boundarylayer}) for scenarios with one edge, one vertex, and two adjacent edges of an element on the boundary. As a result, they yield two, three, three, and four smaller elements, respectively. Note that Templates (b) and (c) introduce a new valence-$3$ EP. In Template (c), if we insert two new vertices on the boundary \cite{liang2009guaranteed}, we obtain Template (d), which could avoid introducing the new valence-3 EP.

\subsection{Supervised Learning}
AFM's nature of using local information not only facilitates the implementation of neural networks for automating and enhancing the rule learning process but also makes it possible to extract datasets from a given mesh. We generate $360$ input boundaries by sampling $N$ points in a unit square domain, where $N \in [4, 100]$ and is an even number, and connect these points sequentially. To avoid invalid geometries, we conduct an intersection check for each newly connected line segment with existing edges. Once we obtain the valid domain boundaries, we use ANSYS to generate corresponding quad meshes for SL training.

The data for SL are input-output pairs extracted at each iteration from the meshes we obtained. In Figure \ref{fig:bigpicture}(b, c), the input tensor of the policy neural network contains $12$ vertices: the reference vertex $P_0$, four vertices on the left $P_i^l, i = 1, 2, 3, 4$, four vertices on the right $P_i^r, i = 1, 2, 3, 4$, and the closest three red vertices $P_i^c, i = 1, 2, 3$, to $P_0$. Note that $P_0, P_i^l, P_i^r$ and $P_i^c$ are all on the front, and the closest three red vertices $P_i^c, i = 1, 2, 3$ are used for intersection checking. We also need valence information of vertices $P_0, P_i^r, P_i^l, i = 1, 2$, because their valences may change during the iteration and turn them from regular vertices into EPs or even cEPs. We define the valence information of these five vertices as ``EP status" to help neural networks distinguish EPs and cEPs from regular vertices. According to the definition of EPs in \cite{wei2022analysis}, we assign the EP status of an interior valence-$V$ vertex to be $V-4$. If the vertex is on the input boundary, its EP status is $-0.5$ when $V<2$ and $V-2$ otherwise. Then we have EP status of $-1$ (interior valence-$3$ EPs), $-0.5$ (regular boundary vertices), $0$ (regular interior or boundary vertices), and $>0$ (EPs). Using these information, the SL and RL neural networks will try to avoid creating EPs when updating the front. There are four updating types for the front, as defined in Figure \ref{fig:accept_reject_types}(a). Type $1$ generates a new vertex, and the number of edges in the front remains the same; Types $2$ and $3$ do not generate new vertices, and they remove two edges from the front; Type $4$ generates two new vertices and inserts two new edges to the front. In the implementation, we reduce the input dimension by normalizing the coordinates. All vertices are transformed by normalizing $P_0$ and $P_1^r$ to $(0, 0)$ and $(1, 0)$ using matrix transformation:

\begin{equation}
\begin{aligned}
    \begin{pmatrix}
         x' \\
         y' \\
         1
    \end{pmatrix}&=
    \begin{pmatrix}
        \cos\theta & \sin\theta & 0 \\
        -\sin\theta & \cos\theta & 0 \\
        0 & 0 & 1
    \end{pmatrix}
    \begin{pmatrix}
        \frac{1}{d} & 0 & 0 \\
        0 & \frac{1}{d} & 0 \\
        0 & 0 & 1
    \end{pmatrix}
    \begin{pmatrix}
        1 & 0 & -x_0 \\
        0 & 1 & -y_0 \\
        0 & 0 & 1
    \end{pmatrix}
    \begin{pmatrix}
        x \\
        y \\
        1
    \end{pmatrix},
\end{aligned}
\end{equation}
where $(x, y), (x', y')$ are the coordinates before and after transformation, $\theta$ is the angle between $P_0P_1^r$ and the horizontal $x-$axis, and $d$ is the edge length $P_0P_1^r$.

\begin{figure}[ht]
\centering
\includegraphics[width=\linewidth]{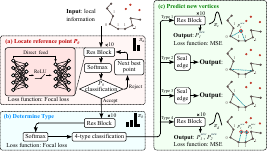}
\vspace{-6mm}
\caption{The residual neural network framework \cite{he2016deep} used in SL. (a) $\pi_a$ is a binary classification neural network that determines the reference vertex $P_0$. (b) $\pi_b$ is a four-class classification neural network that determines the updating type. (c) We then selectively apply $\pi_c$, $\pi_d$ or seal the edge based on $\pi_b$ classification results.}
\label{fig:slarchitecture}
\end{figure}

In Figure \ref{fig:slarchitecture}, we use four residual neural networks to mesh the domain for two reasons: (1) residual connections support very deep neural networks while avoiding the gradient vanishing and over-fitting problem; and (2) a residual block can be easily added to existing neural networks by initializing itself to be an identity mapping, which allows users to add more layers to the neural networks without requiring a complete overhaul of the architecture. The first neural network $\pi_a$ determines whether to select a vertex as the reference $P_0$. All vertices on the front are arranged in ascending order according to the inner angle. We start from small to large and pass local information $S_t$ around the selected vertex into the binary classification neural network $\pi_a$. The softmax layer in $\pi_a$ gives a binary classification result $\left\{\mathscr{P}_{acc}, \mathscr{P}_{rej}\right\}$. The process stops when $\pi_a$ accepts the current vertex to be the reference vertex $P_0$. After we select $P_0$, the local information around $P_0$, $S_t$, is sent to $\pi_b$ to determine the updating type. $\pi_b$ is a four-class classifications neural network. The softmax layer in $\pi_b$ gives a four-class classification result $\left\{\mathscr{P}_1, \mathscr{P}_2, \mathscr{P}_3, \mathscr{P}_4\right\}$ that selects one out of four updating types. Therefore, both $\pi_{a}$ and $\pi_{b}$ adopt Focal Loss \cite{lin2017focal} to address the class imbalance problem. For $\pi_a$ and $\pi_b$, we have
\begin{equation}
\textit{Focal Loss}_{\pi_{a/b}} = -\frac{1}{N_{a/b}}\alpha_{a/b}\sum_{i=1}^{N_{a/b}}\left[\left(1-p^i\right)^\gamma\log p^i\right],
\end{equation}
where $N_{a/b}$ is the total number of neural network $\pi_a$ or $\pi_b$ training rows. For $\pi_a$, $\alpha_{a}$ is a scaling factor of two classes, $p^i$ is the acceptance probability $\mathscr{P}^i_{acc}$ when $\pi_a$ selects the correct category in the ground truth, and the rejection probability $\mathscr{P}^i_{rej}$ otherwise. For $\pi_b$, $\alpha_{b}$ is a scaling factor of four classes, $p^i$ is the acceptance probability $\left\{\mathscr{P}^i_1, \mathscr{P}^i_2, \mathscr{P}^i_3, \mathscr{P}^i_4\right\}$ when the ground truth updating type is $1, 2, 3$ and $4$, respectively. We have $\mathscr{P}^i_{acc}+\mathscr{P}^i_{rej}=1$ and $\mathscr{P}^i_1+\mathscr{P}^i_2+\mathscr{P}^i_3+\mathscr{P}^i_4=1$. $\gamma$ is a tunable parameter that controls the weight of difficult-to-classify samples.

As shown in Figure \ref{fig:slarchitecture}(c), when $\pi_b$ gives Type $1$ classification result, we send local information $S_t$ to the regression neural network $\pi_c$ to generate a new interior vertex $P^{new}_1$. We simply seal the edge when the updating type is $2$ or $3$, and no neural network is needed. When $\pi_b$ gives Type $4$ classification result, we send the same local information $S_t$ to the regression neural network $\pi_d$ to generate two new interior vertices $P^{new}_1$ and $P^{new}_2$. For all $P_1^{new}$ in $\pi_c$ and $P_1^{new}, P_2^{new}$ in $\pi_d$, we utilize their polar coordinates in the following loss function computation. With $P_0$ as the origin, the angle and radius are normalized to $[0, 1]$ by the reference vertex angle $(\angle P_1^rP_0P_1^l)$ and the total length of six line segments around $P_0$ with three on the left and three on the right $(\overline{P_3^rP_2^r}+\overline{P_2^rP_1^r}+\overline{P_1^rP_0}+\overline{P_0P_1^l}+\overline{P_1^lP_2^l}+\overline{P_2^lP_3^l})$. Both $\pi_{c}$ and $\pi_{d}$ adopt MSE loss. For $\pi_c$ we have
\begin{equation}
\textit{MSE}_{\pi_c} = \frac{1}{N_c}\sum_{i=1}^{N_c}\left[\left(\theta_1^i-\hat{\theta}_1^i\right)^2+\left(\rho_1^i-\hat{\rho}_1^i\right)^2\right],
\end{equation}
where $N_c$ is the total number of neural network $\pi_c$ training rows. $\theta_1^i$ and $\hat{\theta}_1^i$ are the prediction and ground truth angle of $P_1^{new}$. $\rho_1^i$ and $\hat{\rho}_1^i$ are the prediction and ground truth radius of $P_1^{new}$. Similarly for $\pi_d$ we have
\begin{equation}
\textit{MSE}_{\pi_d} = \frac{1}{N_d}\sum_{i=1}^{N_d}\left[\left(\theta_1^i-\hat{\theta}_1^i\right)^2+\left(\rho_1^i-\hat{\rho}_1^i\right)^2+\left(\theta_2^i-\hat{\theta}_2^i\right)^2+\left(\rho_2^i-\hat{\rho}_2^i\right)^2\right],
\end{equation}
where $N_d$ is the total number of neural network $\pi_d$ training rows. $\theta_1^i, \hat{\theta}_1^i$ and $\theta_2^i, \hat{\theta}_2^i$ are the prediction and ground truth angle of $P_1^{new}$ and $P_2^{new}$. $\rho_1^i, \hat{\rho}_1^i$ and $\rho_2^i, \hat{\rho}_2^i$ are the prediction and ground truth radii of $P_1^{new}$ and $P_2^{new}$. For all the neural networks, we use $10$ residual blocks \cite{he2016deep} as the middle layers. The size of the training data determines the performance of the model. In principle, the training data should cover the feasible region of the input features comprehensively. We feed 360 input domains or 3.5M training rows to the SL neural networks in the paper.

\subsection{Reinforcement Learning}
To further improve the quad mesh quality generated by the framework, we fine-tune the framework by applying RL to the four trained SL neural networks $\left\{\pi_a, \pi_b, \pi_c, \pi_d\right\}$ using five randomly selected new boundaries. The fine-tuned RL neural networks $\left\{\pi'_a, \pi'_b, \pi'_c, \pi'_d\right\}$ have the same architecture as Figure \ref{fig:slarchitecture}. The whole RL algorithm is demonstrated in Algorithm \ref{alg:RL}.

\begin{algorithm}
    \textbf{Input}: SL neural networks $\left\{\pi_a, \pi_b, \pi_c, \pi_d\right\}$, and 5 randomly selected boundaries $B^0_i$\\
    \textbf{Output}: RL neural networks $\left\{\pi'_a, \pi'_b, \pi'_c, \pi'_d\right\}$ \\
    \textbf{Goal}: Improve the neural network meshing performance
    \caption{Reinforcement Learning}
    \label{alg:RL}
    \begin{algorithmic}[1]
        \State Copy the SL neural networks and denote them as initial RL neural networks $\left\{\pi'_a, \pi'_b, \pi'_c, \pi'_d\right\}$
        \For{episode $i = 1, 2, \cdots, M$}
        \State Get an initial state $S_1^{ij}$ ($i$ is omitted hereafter)
        \For{time step $t = 1, 2, \cdots, M$}
        \State $A_t\leftarrow$ sampling action $(\pi'_a, \pi'_b, \pi'_c, \pi'_d, S_t, \epsilon)$
        \State Form a new quad element and update the front boundary with $A_t$
        \State Get the reward (element quality) $R_t=R^s_tR^{ep}_t$ and the next state $S_{t+1}$
        \If {the number of edges on the remaining front boundary $=4$}
        \State Break
        \EndIf
        \EndFor
        \If {$R^{fin} = \frac{1}{M} \sum_{i=1}^{M} R_i + \min\{R_1, R_2, \cdots, R_M\} > R'_{fin}$}
        \State Extract dataset from the new RL mesh
        \State Update $\left\{\pi'_a, \pi'_b, \pi'_c, \pi'_d\right\}$ weights with the new dataset
        \EndIf
        \EndFor
    \end{algorithmic}
\end{algorithm}

To obtain meshes with higher quality than the SL-generated meshes, we introduce additional exploration to the neural network-guided action by adding noise $\epsilon$ to RL neural networks. We add Dirichlet noise $Dir(\alpha)$, a type of probability distribution that assigns probabilities to an arbitrary number of outcomes, to classification neural networks $\pi'_a$ and $\pi'_b$ and $2-D$ Gaussian noise $\mathcal{N}(\mu, \Sigma)$ to regression neural networks $\pi'_c$ and $\pi'_d$. As a supplement to the fifth line of Algorithm \ref{alg:RL}, we achieve this by sampling four actions from four different distributions. We have
\begin{align}
A^{a/b}_t&\sim\eta \pi'_{a/b}(S_t) + (1-\eta)Dir(\alpha^{a/b}),\\
A^{c/d}_t&\sim\mathcal{N}\left(\pi'_{c/d}(S_t), \Sigma^{c/d}\right).
\end{align}
Based on several trials, we set $\eta=0.99, \alpha^{a/b}=0.05,$ and $\Sigma^{c/d}=0.0001I_{2\times2}$ to balance the neural network exploration and exploitation trade-off. Here $I_{2\times2}$ is a $2\times2$ identity matrix.

The RL algorithm takes SL policy neural networks and five randomly selected initial boundaries as input and outputs fine-tuned RL policy neural networks. The neural network $\pi'_a$, at time step $t$, selects a vertex as the reference $P_0$ and records the state $S_t$ describing vertices near $P_0$ from the environment. Then $\pi'_b, \pi'_c\ \text{or}\ \pi'_d$ conduct an action $A(S_t)$ that determines the updating type and assigns new vertices (if any) to the environment. The environment responds to the action and transits into a new state  $S_{t+1}$ at time step $t+1$. After meshing, we do not apply any post-processing. Finally, we calculate the reward related to the mesh quality. If the reward exceeds the previous mesh, we save the current mesh and extract the dataset.

The reward function is critical in guiding the mesh optimization direction in RL. In our framework, we have two reward functions: the squareness reward function $R^s$ and the EP penalty reward function $R^{ep}$. $R^s$ is defined based on inner angles $\theta_1, \theta_2, \theta_3, \theta_4$ and edge lengths $l_1, l_2, l_3, l_4$ of a quad. We have
\begin{equation}
\begin{aligned}
    R^s = \sqrt[3]{\frac{\min\{\theta_1, \theta_2, \theta_3, \theta_4\}}{\ang{90}}\left(2 - \frac{\max\{\theta_1, \theta_2, \theta_3, \theta_4\}}{\ang{90}}\right)\frac{\min\{l_1, l_2, l_3, l_4\}}{\max\{l_1, l_2, l_3, l_4\}}}.
\end{aligned}
\end{equation} $R^s$ is bounded between $0$ and $1$. The higher the $R^s$, the more similar a quad is to a square, regardless of the element size. 

Minimizing the number of EPs and cEPs is also important for high-order finite element analysis such as IGA \cite{wei2022analysis}. We assign reward ``1'' to all regular vertices and penalty ``0'' to all EPs and cEPs as shown in Figure \ref{fig:bigpicture}(c), and then compute the number of EPs ($N_{ep}$) and close EP pairs ($N_{cep}$) in a mesh with a total of $N_{tot}$ vertices. $R^{ep}$ is defined as

\begin{equation}
\begin{aligned}
    R^{ep} = \left(1 - \frac{N_{ep}}{N_{tot}}\right)\left(1 - \frac{N_{cep}}{N_{tot}}\right).
\end{aligned}
\end{equation} With both $R^s$ and $R^{ep}$ defined, the total reward is evaluated as the sum of the mean and the minimum reward value across all generated quads:

\begin{equation}\label{eq:fin}
\begin{aligned}
    R^{fin} = \frac{1}{M} \sum_{i=1}^{M} R^s_iR^{ep}_i + \min\{R^s_1R^{ep}_1, R^s_2R^{ep}_2, \cdots, R^s_MR^{ep}_M\},
\end{aligned}
\end{equation} where $M$ is the number of quads produced.

The most significant advantage of RL in the SRL-assisted AFM is that we can arbitrarily define the reward function. In this research, we use Equation \ref{eq:fin} to optimize the number and distribution of EPs and quad squareness simultaneously. As training processes, the mesh quality improves by maximizing the reward function, and finally the quality will significantly exceeds previous performance.
\section{Numerical Results and Discussion}
\begin{figure}[ht]
\centering
\includegraphics[width=\linewidth]{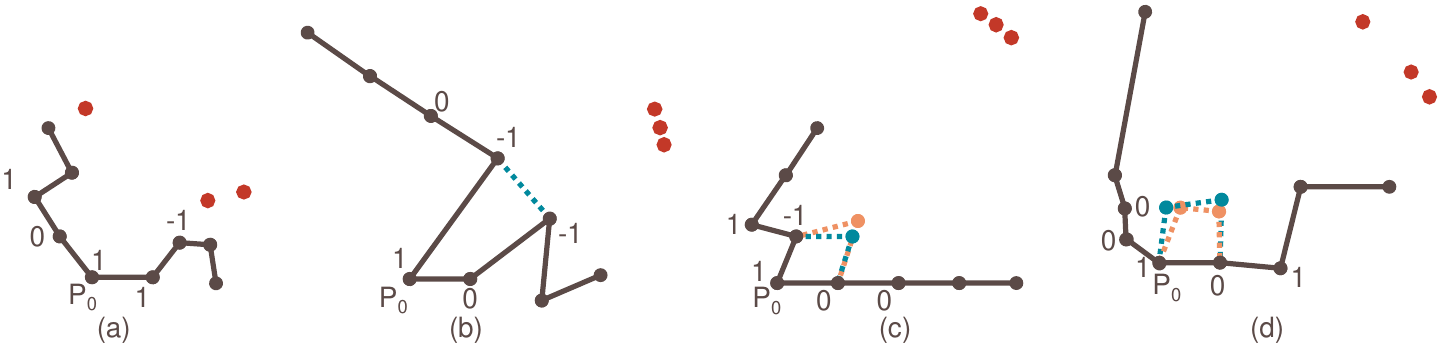}
\vspace{-6mm}
\caption{SL policy neural network outputs of four held-out examples in Table \ref{tab:SLnn}. The local vertices $P_0, P_i^r, P_i^l, i=1,2,3,4$, are in black. Three close vertices $P_i^c, i = 1,2,3$, are in red. EP status are listed next to the corresponding vertices. The new vertices predicted by the neural network (if any) are in blue. The new vertices from ground truth (if any) are in orange. (a) $\pi_a$ accepts $P_0$ as the reference vertex. (b) $\pi_b$ selects updating type $2$ and two existing points are connected to seal the edge. (c) $\pi_b$ selects updating type $1$, a new point $P_1^{new}$ is generated based on $\pi_c$ prediction. (d) $\pi_b$ selects updating type $4$, two new points $P_1^{new}, P_2^{new}$ are generated based on $\pi_d$ predictions.}
\label{fig:sl}
\end{figure}

The meshes produced by ANSYS serve as the basis for extracting the SL training dataset. In the SL phase, we train four SL policy neural networks $\left\{\pi_a, \pi_b, \pi_c, \pi_d\right\}$ using approximately $3.5M$ training rows. Figure \ref{fig:sl} shows SL policy neural network outputs of four held-out examples extracted from the training dataset. With $10$ residual blocks, the SL neural networks achieve accuracies of $99.78\%, 98.07\%$ in choosing the reference $P_0$ and front updating types, and MSE of $0.0017, 0.0037$ in predicting Type $1$ and Type $4$ new vertices, respectively (Table \ref{tab:SLnn}). In these testing cases, $\pi_a$ and $\pi_b$ give the correct classification with probabilities of $91.0\%$ and $98.0\%$. $\pi_c$ and $\pi_d$ also give good predictions of new vertices; see the blue and orange points in Figure \ref{fig:sl}(c,d). In the RL phase, to generate enough high quality datasets for RL training, we apply the RL Algorithm \ref{alg:RL} to five randomly selected complex domains from the $360$ input training boundaries in Figure \ref{fig:bigpicture}(a). The training starts from SL neural networks $\left\{\pi_a, \pi_b, \pi_c, \pi_d\right\}$ and each domain is meshed with many episodes in $24$ hours to collect high-quality meshes. Over the course of RL training, more than $10,000$ meshes are generated for each domain, in which $9\%$ $(1,130)$ meshes have their quality exceeding the corresponding SL mesh quality. We keep $500$ meshes with the highest mesh quality $R^{fin}$, and extract dataset from these $500$ meshes to train the RL policy neural networks $\left\{\pi'_a, \pi'_b, \pi'_c, \pi'_d\right\}$. In both the SL and RL phases, we use stochastic gradient descent as the gradient optimizer \cite{bottou2012stochastic}, and the learning rate is fixed at $10^{-5}$.

\vspace{-5mm}
\begin{table}[ht]
\centering
\resizebox{\linewidth}{!}{\begin{threeparttable}
\caption{Examples and statistics of SL neural network outputs with $N=10$ residual blocks.}
\label{tab:SLnn}
\begin{tabular}{c|cccc}
\toprule
Neural Network & $\pi_a(\%)$ & $\pi_b(\%)$ & $\pi_c$ & $\pi_d$\\
Raw Output& $(\mathscr{P}_{acc}, \mathscr{P}_{rej})$ & $(\mathscr{P}_1, \mathscr{P}_2, \mathscr{P}_3, \mathscr{P}_4)$ & $(\theta_1, \rho_1)$& $(\theta_1, \rho_1), (\theta_2, \rho_2)$\\\hline
Prediction & $(91.0, 9.0)$ & $(1.9, 98.0, 0.0, 0.1)$ & $(0.46, 0.25)$ & $(0.57, 0.18), (0.32, 0.27)$\\
Ground Truth & $(100.0, 0.0)$ & $(0.0, 100.0, 0.0, 0.0)$ & $(0.55, 0.30)$ & $(0.48, 0.19), (0.28, 0.24)$\\\hline
Accuracy/MSE & $99.78$ & $98.07$ & $0.0017$ & $0.0037$\\
\bottomrule
\end{tabular}
\end{threeparttable}}
\end{table}

\begin{figure}[ht]
\centering
\includegraphics[width=\linewidth]{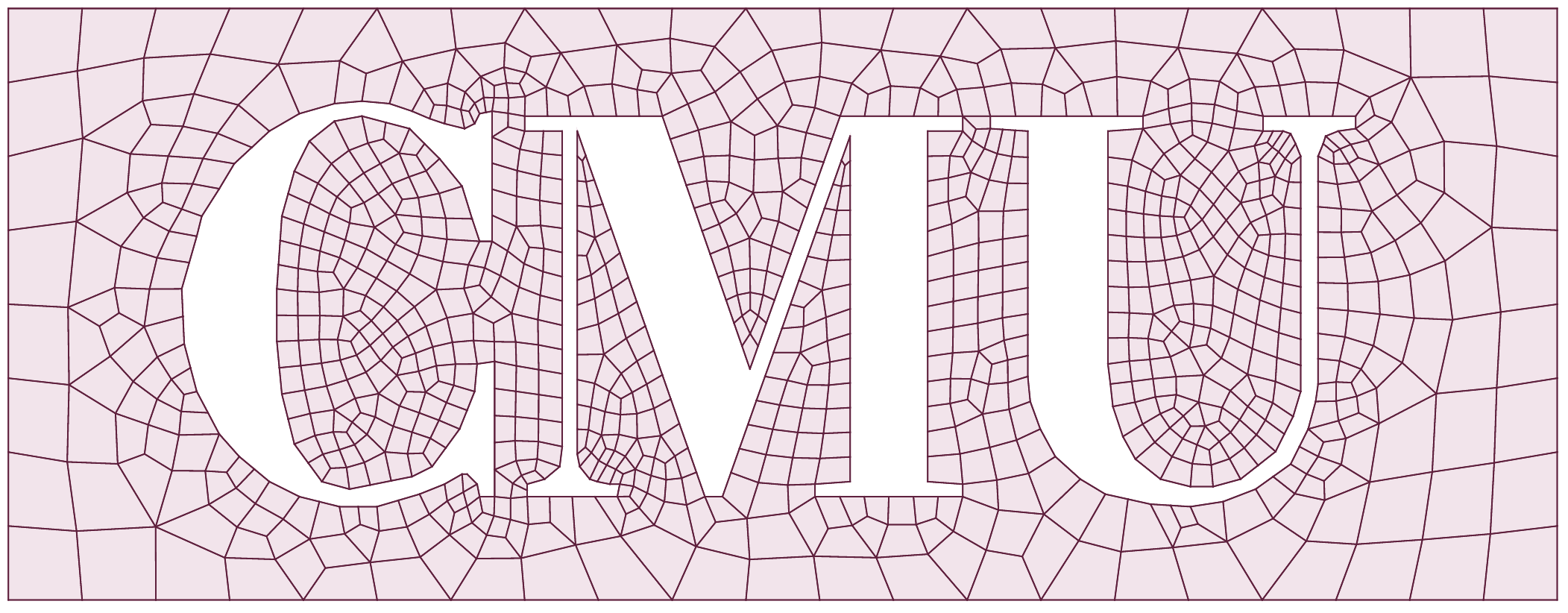}
\vspace{-6mm}
\caption{The CMU logo mesh covering the exterior domain is formed by stripping off the logo in a rectangular box.}
\label{fig:cmusl}
\end{figure}

\begin{figure}[ht]
\centering
\includegraphics[width=\linewidth]{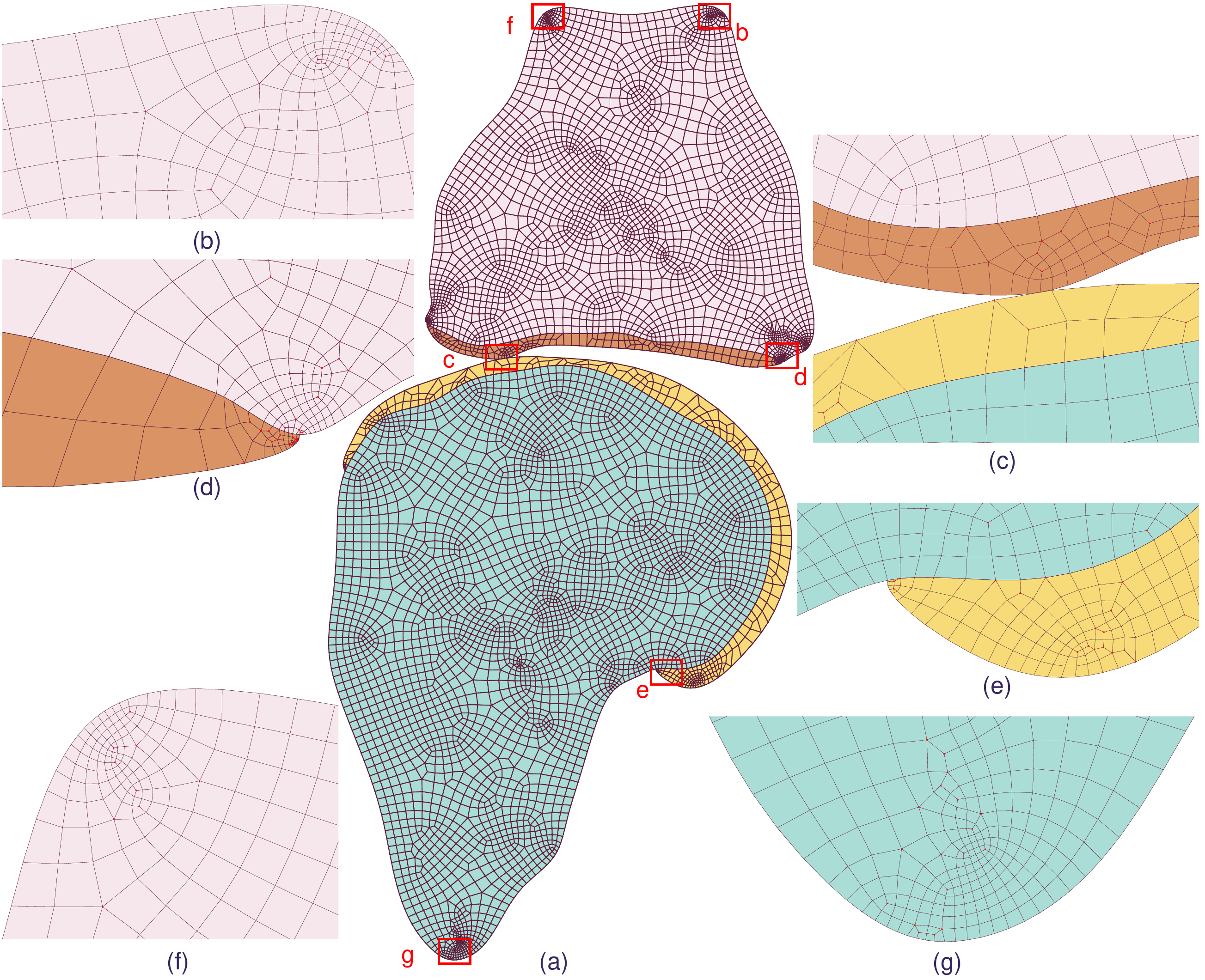}
\vspace{-6mm}
\caption{The knee joint cross-section. (a) Multiple meshes representing the femur (pink), femoral cartilage (orange), tibial cartilage (yellow), and tibia (cyan). (b-f) Zoom-in pictures of red boxes in (a).}
\label{fig:kneebotsl}
\end{figure}

\begin{figure}[ht]
\centering
\includegraphics[width=\linewidth]{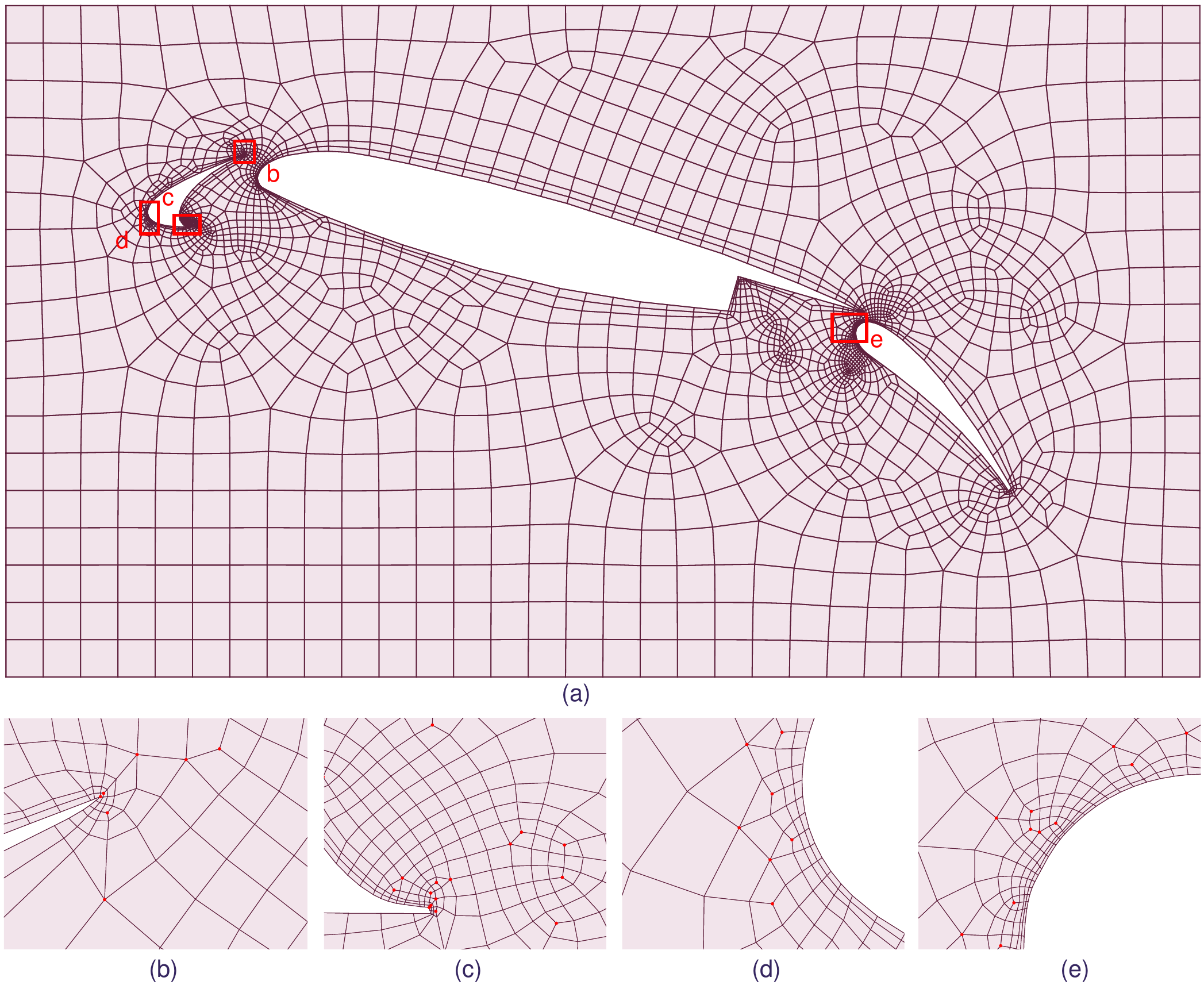}
\vspace{-6mm}
\caption{The airfoil with three components. (a) The mesh covering the exterior domain is formed by stripping off the airfoil in a rectangular box. Two boundary layers are added to the airfoil surfaces. (b-e) Zoom-in pictures of (a).}
\label{fig:airfoilsl}
\end{figure}

\begin{figure}[ht]
\centering
\includegraphics[width=\linewidth]{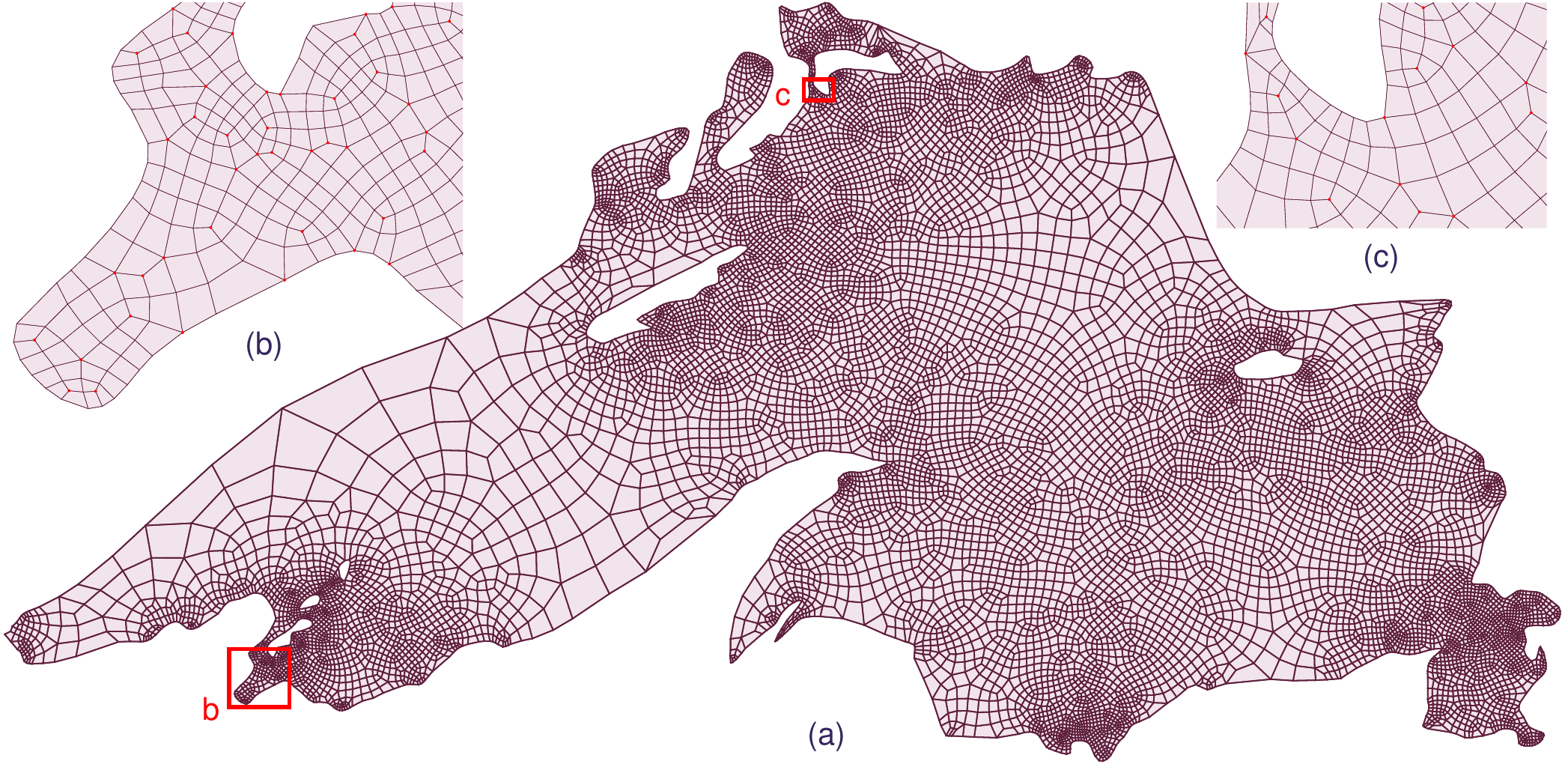}
\vspace{-6mm}
\caption{The Lake Superior map with multiple holes and unbalanced seeds. (a) Final all-quad mesh. (b-c) Zoom-in pictures of the red boxes in (a).}
\label{fig:lakesl}
\end{figure}

We apply our SRL-assisted AFM framework to five complicated domains. Each boundary domain inherits some unique challenges. The first curve is shown in Figure \ref{fig:bigpicture}(a) with SL and RL meshes. Figure \ref{fig:cmusl} shows the CMU logo with sharp angles and narrow regions. Our framework meshes the domain successfully within $1$ second. In Figure \ref{fig:kneebotsl}, we take one slice of the segmented mask from a knee MRI dataset \cite{ambellan2019automated} and construct four knee joint components (femur, femoral cartilage, tibial cartilage, and tibia) from the mask. The constructed mesh conforms to each other exactly on shared boundaries while exhibiting good adaptivity. The mesh generation for this example takes around $10$ seconds, slightly slower than meshing a whole boundary due to the time required to assign conformal seeds. In Figure \ref{fig:airfoilsl}, we apply our framework to an airfoil consisting of three sharp-angle components. It takes $8$ seconds to generate an adaptive mesh that effectively preserves sharp angles and narrow regions. Our framework exhibits remarkable efficacy in handling inputs at a large scale while simultaneously reducing the number of elements critical for simulations. In addition, our model supports adding arbitrary numbers of boundary layers. Here, we add two boundary layers to the airfoil boundaries to illustrate the effectiveness of this feature. The last and most complex mesh is the Lake Superior map in Figure \ref{fig:lakesl}. To the best of our knowledge, no ML-based mesh generation method has been attempted on such complex boundaries with multiple holes, sharp angles, narrow regions, and unbalanced seeds before. Our framework managed to mesh the domain in $118$ seconds, showing superior stability in large-scale meshes. All the final meshes preserve correct topology structures, and no intersection detection or post-processing optimization are used when generating these meshes.

Table \ref{tab:quality} shows statistics of our resulting meshes. In SRL-assisted AFM, the mesh size is proportional to the seed density. Our mesh size has the same scale as the previous guarantee-quality methods \cite{liang2009guaranteed, liang2011hexagon, liang2012matching}. Among the five testing meshes, the proportions of the adopted four updating types $f_1, f_2, f_3, f_4$ are roughly equivalent to $90\%, 4.5\%, 4.5\%, 1\%$, which means $f_1\gg f_2\approx f_3\gg f_4$. In addition, the numbers of EPs and close EPs are also very stable, around $11\%$ and $7\%$. In SL meshes, the percentages are around $17\%$ and $12\%$. Instead of guaranteeing the minimum and maximum angles, the AFM directly uses neural network planning, resulting in slightly poor mesh quality in terms of angle range and Jacobian metric on complex domains with sharp features. Additionally, the mesh has superior aspect ratio performance compared to guarantee-quality methods due to adaptive seeding in pre-processing period and new vertex prediction by $\pi'_c, \pi'_d$. Finally, the meshing time is proportional to the geometry complexity. The Lake Superior example in Figure \ref{fig:lakesl} is the most complex example because it has massive amount of unbalanced seeds, narrow regions, and sharp angles. Time records reveal that in Algorithm \ref{alg:AFM}, determining the reference vertex (Steps $12- 15$) is the most time-consuming step since it iterates all vertices on the evolving front to find the proper reference vertex. All the results were computed on a PC with an Intel i$7-12700$ CPU and $64$GB memory. The code is written in Python and available at \url{https://github.com/CMU-CBML/SRL-AssistedAFM}.

\begin{table}[ht]
\centering
\resizebox{\linewidth}{!}{\begin{threeparttable}
\caption{Statistics of the resulting meshes.}
\label{tab:quality}
\begin{tabular}{c|cccccc}
\toprule
Domain & Mesh Size & Aspect Ratio & Valence& Angle & Jacobian & Time\\
& [Vert$\#$, Elem$\#$] & [Best, Worst] & [EP, cEP$^*$] & [Min, Max] & [Worst, Best] & (s)\\ \hline
Curve & $[1,401, 1,520]$ & $[1.0, 3.9]$ & $[168, 98]$ & $[35^\circ, 148^\circ]$ & $[0.75, 1.0]$ & $0.4$\\
CMU Logo & $[924, 1,115]$ & $[1.0, 3.5]$ & $[111, 92]$ & $[24^\circ, 140^\circ]$ & $[0.68, 1.0]$ & $1.1$\\
Knee Joint & $[6,196, 6,617]$ & $[1.0, 2.4]$ & $[620, 372]$ & $[43^\circ, 137^\circ]$ & $[0.72, 1.0]$ & $10.9$\\
Air Foil & $[2,782, 2,953]$ & $\ \ [1.0, 4.8]^+$ & $[250, 167]$ & $[29^\circ, 150^\circ]$ & $[0.68, 1.0]$ & $8.1$\\
Lake Superior & $[12,150, 11,618]$ & $[1.0, 7.2]$ & $[389, 223]$ & $[15^\circ, 156^\circ]$ & $[0.60, 1.0]$ & $118.1$\\
\bottomrule
\end{tabular}
\begin{tablenotes}
      \item[*] cEP is the number of pairs of EPs that are adjacent to each other.
      \item[+] After adding boundary layers, the worst aspect ratio becomes $16.0$.
\end{tablenotes}
\end{threeparttable}}
\end{table}

\section{Conclusion and Future Work}
This paper presents a novel method that automatically generate quad meshes for complex planar domains with four neural networks. The computational framework SRL-assisted AFM integrates the SL-RL algorithm with the AFM. 
\begin{itemize}
    \item We generate a large number of planar domain meshes ($3.5M$ training rows) to increase the diversity of the training dataset. As a result, our framework can mesh new boundaries that are far more complex than those tested in existing ML-based mesh generation methods.
    \item We do not adopt any quality improvement or intersection detection module throughout the pipeline, indicating that our framework can correct errors. This capability significantly simplifies the AFM pipeline and improves the meshing efficiency. For the tested five domains, SRL-assisted AFM generates $\sim\!\!1,000$ quads per second in average.
    \item Our SRL-assisted AFM achieves high mesh quality that rivals other methods without the need of post-processing operations or prior knowledge. The framework can generate adaptive meshes to reduce computation burden. The angle range and scaled Jacobian metrics are close to previous guarantee-quality methods. It also significantly improves the aspect ratio and EP penalty metric.
    \item We define the reward function to reduce the number of EPs and adjacent EPs. The mesh adaptivity is achieved by a size function that assigns seeds according to local front boundary curvature. Users can also assign boundary layers after the mesh is generated.
\end{itemize}

In the current implementation, the code is written in Python. As a result, it is relatively slower than C$/$C++. In the future, we will include all modules in C++ to improve computational efficiency. We also aim to extend this SRL-assisted AFM framework to curved surface mesh generation and other meshing approaches beyond the AFM. We anticipate that this technology will significantly improve and boost data-driven mesh generation.

\section{Acknowledgment}
H. Tong, K. Qian, and Y. J. Zhang were supported in part by the NSF grant CMMI-1953323 and a Honda grant. K. Qian was also supported by Bradford and Diane Smith Graduate Fellowship. This work used RM-node and GPU-node on Bridges-2 Supercomputer at Pittsburgh Supercomputer Center \cite{xsede,ecss} through allocation ID eng170006p from the Advanced Cyberinfrastructure Coordination Ecosystem: Services \& Support (ACCESS) program, which is supported by NSF grants \#2138259, \#2138286, \#2138307, \#2137603, and \#2138296.

\vspace{-3mm}
\bibliographystyle{elsarticle-num}
\bibliography{refs}
\end{document}